\documentclass{article}
\usepackage{amssymb}
\usepackage{amsfonts}
\usepackage{amsmath}
\usepackage{graphicx}

\input{tcilatex}

\begin{document}

\title{Runge-Kutta Methods: Local error control does not imply global error
control}
\author{J.S.C. Prentice \\
Department of Applied Mathematics, \\
University of Johannesburg, \\
South Africa}
\maketitle

\begin{abstract}
We study the relationship between local and global error in Runge-Kutta
methods for initial-value problems in ordinary differential equations. We
show that local error control by means of local extrapolation does not
equate to global error control. Our analysis shows that the global error of
the higher-order solution is propagated under iteration, and this can cause
an uncontrolled increase in the global error of the lower-order solution. We
find conditions under which global error control occurs during the initial
stages of the RK integration, but even in such a case the global error is
likely to eventually exceed the user-defined tolerance.
\end{abstract}

\textbf{Keywords:} Runge-Kutta; Local error; Global error; Local
extrapolation; Error control

\textbf{MSC 2010:} 65L05, 65L06, 65G99

\section{Introduction}

Runge-Kutta (RK) methods are the most popular choice of one-step methods for
solving problems of the form%
\begin{eqnarray}
y^{\prime } &=&f\left( x,y\right)  \label{ivp} \\
y\left( x_{0}\right) &=&y_{0}  \notag
\end{eqnarray}%
numerically. In the implementation of these methods, local error control via
local extrapolation is the preferred choice of error control. It is known,
however, that this form of error control does not amount to control of the
global error. In this paper, we seek to give our own interpretation and
discussion of this problem, including an analysis showing that global error
control, if it does occur, is limited under RK iteration.

\section{Relevant concepts}

We now define local and global errors in RK methods formally, and study the
propagation of local error in the implementation of an explicit RK method.
We will also find the relationship between local and global error.

\subsection{Runge-Kutta methods}

The most general definition of a Runge-Kutta method \cite{RKdefn} is%
\begin{equation}
\begin{array}{l}
k_{p}=f\left( x_{i}+c_{p}h,w_{i}+h\sum\limits_{q=1}^{m}a_{pq}k_{q}\right) 
\text{ \ \ \ \ \ \ \ }p=1,2,...,m \\ 
w_{i+1}=w_{i}+h\sum\limits_{p=1}^{m}b_{p}k_{p}%
\end{array}
\label{RK defn}
\end{equation}%
Such a method is said to have $m$ stages (the $k_{q}$). We note that if $%
a_{pq}=0$ for all $p\leqslant q,$ then the method is said to be \textit{%
explicit}; otherwise, it is known as an \textit{implicit} RK method. The
number of stages is related to the order of the method. The symbol $w$ is
used here and throughout to indicate the approximate numerical solution,
whereas the symbol $y$ will denote the true solution.

We denote an explicit Runge-Kutta method of order $z$ (RK$z)$ for solving (%
\ref{ivp}) by%
\begin{equation*}
w_{i+1}^{z}=w_{i}^{z}+hF^{z}\left( x_{i},w_{i}^{z}\right)
\end{equation*}%
where $w_{i}^{z}$ denotes the numerical approximation to $y\left(
x_{i}\right) $ and $F^{z}\left( x,y\right) $ is a function associated with
the particular RK$z$ method. Indeed, $F^{z}$ is simply the linear
combination of the stages for that particualr method, as in%
\begin{equation*}
F^{z}=\sum\limits_{p=1}^{m}b_{p}k_{p}.
\end{equation*}

\subsection{Local and global errors}

We define the global error in a numerical solution at $x_{i+1}$ by

\begin{equation*}
\Delta _{i+1}^{z}\equiv w_{i+1}^{z}-y_{i+1},
\end{equation*}%
and the local error at $x_{i+1}$ by 
\begin{equation}
\varepsilon _{i+1}^{z}\equiv \left[ y_{i}+hF^{z}\left( x_{i},y_{i}\right) %
\right] -y_{i+1}.  \label{eps local}
\end{equation}%
In the above, $y_{i}$ denotes the true solution $y\left( x_{i}\right) ,$ and
similarly for $y_{i+1}.$ Note the use of the exact value $y_{i}$ in the
bracketed term in (\ref{eps local}).

Previously, we have shown \cite{jscp1} that%
\begin{eqnarray}
\Delta _{i+1}^{z} &=&\varepsilon _{i+1}^{z}+\alpha _{i}^{z}\Delta _{i}^{z}
\label{master error equn} \\
\alpha _{i}^{z} &\equiv &1+hF_{y}^{z}\left( x_{i},\xi _{i};h\right) ,  \notag
\end{eqnarray}%
where $\xi _{i}\in \left( y_{i},y_{i}+\Delta _{i}^{z}\right) $. Equation (%
\ref{master error equn}) provides the relationship between local and global
errors in RK$z$. We will assume that $\Delta _{0}=0$ (i.e. the initial value
is known exactly). We see that the global error at any node $x_{i+1}$ is the
sum of a local error term and a term incorporating the global error at the
previous node. For convenience, we will drop the argument from $%
F_{y}^{z}\left( x_{i},\xi _{i};h\right) $ in the remainder of the paper; its
presence is implied.

For RK$z$, it is well-known that%
\begin{eqnarray*}
\varepsilon _{i+1}^{z} &\propto &h^{z+1} \\
\Delta _{i+1}^{z} &\propto &h^{z}.
\end{eqnarray*}

\subsection{Local error estimation}

Consider two RK methods of order $z$ and $z+1$. Let $w_{i+1}^{z}$ denote the
approximate solution at $x_{i+1}$ obtained with the order $z$ method, and
similarly for $w_{i+1}^{z+1}.$ Let the local error at $x_{i+1}$ in the order 
$z$ method be denoted by $\varepsilon _{i+1}^{z}=\beta _{i+1}^{z}h^{z+1},$
and similarly for $\varepsilon _{i+1}^{z+1}=\beta _{i+1}^{z+1}h^{z+2}.$
Hence, with $w_{i}^{z},w_{i}^{z+1}=y_{i},$ we have%
\begin{eqnarray*}
w_{i+1}^{z}-w_{i+1}^{z+1}=\varepsilon _{i+1}^{z}-\varepsilon _{i+1}^{z+1}
&=&\beta _{i+1}^{z}h^{z+1}-\beta _{i+1}^{z+1}h^{z+2} \\
&\approx &\beta _{i+1}^{z}h^{z+1}
\end{eqnarray*}%
if $h$ is sufficiently small. This gives%
\begin{equation}
\beta _{i+1}^{z}\approx \frac{w_{i+1}^{z}-w_{i+1}^{z+1}}{h^{z+1}}.
\label{beta local extrap}
\end{equation}

\subsection{Local error control}

Once we have estimated the local error, we can perform error control. Assume
that we require that the local error at each step must be less than a
user-defined tolerance $\delta .$ Moreover, assume that, using stepsize $h$,
we find%
\begin{equation*}
\left\vert \varepsilon _{i+1}^{z}\right\vert =\left\vert \beta
_{i+1}^{z}h^{z+1}\right\vert \geqslant \delta .
\end{equation*}%
In other words, the magnitude of the local error $\varepsilon _{i+1}^{z}$
exceeds the desired tolerance. We remedy the situation by determining a new
stepsize $h^{\ast }$ from%
\begin{equation}
\left\vert \beta _{i+1}^{z}\left[ h^{\ast }\right] ^{z+1}\right\vert =\delta
\Rightarrow h^{\ast }=\left( \frac{\delta }{\left\vert \beta
_{i+1}^{z}\right\vert }\right) ^{\frac{1}{z+1}}  \label{h*=...}
\end{equation}%
and we repeat the RK computation with this new stepsize. This, of course,
gives%
\begin{equation*}
x_{i+1}=x_{i}+h^{\ast }.
\end{equation*}%
This procedure is then carried out on the next step, and so on. Such form of
error control is known as \textit{absolute} error control. If the estimated
error does not exceed the tolerance, then no stepsize adjustment is
necessary, and we proceed to the next step.

Often, we introduce a so-called `safety factor' $\sigma ,$ as in%
\begin{equation*}
h^{\ast }=\sigma \left( \frac{\delta }{\left\vert \beta
_{i+1}^{z}\right\vert }\right) ^{\frac{1}{z+1}}
\end{equation*}%
where $\sigma <1,$ so that the new stepsize is slightly smaller than that
given by (\ref{h*=...}). This is an attempt to cater for the possibility
that $\beta _{i+1}^{z}$ may have been underestimated, due to the assumptions
made in deriving (\ref{beta local extrap}). The choice of the value of $%
\sigma $ is subjective, although a representative value is $0.8.$

Note that because this error control algorithm is applied on each step, we
could find that over the interval of integration we have stepsizes of
varying lengths. For this reason, it is appropriate to make the replacement%
\begin{equation*}
h\rightarrow h_{i}\equiv x_{i+1}-x_{i}
\end{equation*}%
in (\ref{RK defn}).

\subsection{Propagation of the higher-order solution}

There is a very important point that must be discussed. The method for
determining $\beta _{i+1}^{z}$ hinged on the requirement $%
w_{i}^{z},w_{i}^{z+1}=y_{i}.$ However, we only have the exact solution at
the initial point $x_{0};$ at all subsequent nodes, the solution is
approximate. How do we meet the requirement $w_{i}^{z},w_{i}^{z+1}=y_{i}?$

In the case of local extrapolation, the answer is simple: simply use the
higher-order solution $w_{i}^{z+1}$ as input to generate both $w_{i+1}^{z}$
(with the order $z$ method), and $w_{i+1}^{z+1}$ (with the order $z+1$
method). In other words, we are assuming that $w_{i}^{z+1}$ is accurate
enough, relative to $w_{i}^{z},$ to be regarded as the exact value, an
assumption entirely consistent with the assumption made in deriving (\ref%
{beta local extrap}). This means that we determine the higher-order solution
at each node, and this solution is used as input for both methods in
computing solutions at the next node. The question of whether or not the
global error that accumulates in the higher-order solution affects the
calculation of $\beta _{i+1}^{z}$ in (\ref{beta local extrap}) is addressed
in the next section.

\section{Analysis}

\subsection{The problem}

Now, as per the last paragraph of the preceding section, assume that $%
w_{i}^{z+1}$ is used to generate $w_{i+1}^{z}$ and $w_{i+1}^{z+1}.$ Such
value of $w_{i+1}^{z}$ (and associated quantities) will be denoted $%
w_{i+1}^{z,z+1}$. Hence, we have%
\begin{eqnarray*}
\Delta _{i+1}^{z,z+1} &=&\beta _{i+1}^{z}h_{i}^{z+1}+\alpha
_{i}^{z,z+1}\Delta _{i}^{z+1} \\
&=&\beta _{i+1}^{z}h_{i}^{z+1}+\Delta _{i}^{z+1}+h_{i}F_{y}^{z,z+1}\Delta
_{i}^{z+1} \\
\Delta _{i+1}^{z+1} &=&\beta _{i+1}^{z+1}h_{i}^{z+2}+\alpha _{i}^{z+1}\Delta
_{i}^{z+1} \\
&=&\beta _{i+1}^{z+1}h_{i}^{z+2}+\Delta _{i}^{z+1}+h_{i}F_{y}^{z+1}\Delta
_{i}^{z+1}.
\end{eqnarray*}

Thus,%
\begin{eqnarray}
w_{i+1}^{z,z+1}-w_{i+1}^{z+1} &=&\beta _{i+1}^{z}h_{i}^{z+1}+\Delta
_{i}^{z+1}+h_{i}F_{y}^{z,z+1}\Delta _{i}^{z+1}  \notag \\
&&-\left( \beta _{i+1}^{z+1}h_{i}^{z+2}+\Delta
_{i}^{z+1}+h_{i}F_{y}^{z+1}\Delta _{i}^{z+1}\right)  \notag \\
&=&\beta _{i+1}^{z}h_{i}^{z+1}-\beta _{i+1}^{z+1}h_{i}^{z+2}+\left(
F_{y}^{z,z+1}-F_{y}^{z+1}\right) h_{i}\Delta _{i}^{z+1}
\label{cancellation of GE} \\
&\approx &\beta _{i+1}^{z}h_{i}^{z+1}  \notag
\end{eqnarray}%
for small $h_{i}$, because $h_{i}\Delta _{i}^{z+1}=O\left(
h_{i}^{z+2}\right) .$ We see that the presence of global error in the
higher-order solution does not affect the expression for $\beta _{i+1}^{z}$
obtained under the assumption $w_{i}^{z},w_{i}^{z+1}=y_{i},$ particularly if 
$h_{i}$ is small.

However, the expression for $\Delta _{i+1}^{z,z+1}$ informs of a potential
problem: we have%
\begin{equation}
\Delta _{i+1}^{z,z+1}=\beta _{i+1}^{z}h_{i}^{z+1}+\alpha _{i}^{z,z+1}\Delta
_{i}^{z+1},  \label{del z z+1}
\end{equation}%
where $\Delta _{i}^{z+1}$ is the global error in $w_{i}^{z+1}.$ In (\ref%
{cancellation of GE}), we see that a subtractive cancellation ensures that
the $\Delta _{i}^{z+1}$ term does not enter directly into the estimate for $%
\beta _{i+1}^{z}.$ Nevertheless, even if $\left\vert \beta
_{i+1}^{z}h_{i}^{z+1}\right\vert \leqslant \delta ,$ we could still have $%
\left\vert \Delta _{i+1}^{z,z+1}\right\vert >\delta ,$ perhaps substantially
so, if $\left\vert \Delta _{i}^{z+1}\right\vert $ is large. Moreover, we
should certainly expect that $\left\vert \Delta _{i}^{z+1}\right\vert $
could become large under iteration (i.e. as $i$ increases), since global
error is essentially an accumulation of local errors. The point here is
that, even if local error control is effective, the global error $\Delta
_{i+1}^{z,z+1}$ could become large, and could grow in an uncontrolled
fashion.

\subsection{Bounded global error via local error control?}

Let us investigate the effect on the global error if local error control,
via local extrapolation, is implemented. Consider the expression obtained
previously for the global error at $x_{i+1}$%
\begin{equation}
\Delta _{i+1}^{z}=\varepsilon _{i+1}^{z}+\alpha _{i}^{z}\Delta _{i}^{z}.
\label{del n = eps + alpha del}
\end{equation}%
We assume $\Delta _{0}=0.$ If we have the exact value $y_{i}$ at each node,
then we have%
\begin{equation*}
\Delta _{i+1}^{z}=\varepsilon _{i+1}^{z}
\end{equation*}%
at each node, so that the global error is equal to the local error. If the
local error has been controlled (subject to tolerance $\delta ),$ we have%
\begin{equation*}
\left\vert \Delta _{i+1}^{z}\right\vert \leqslant \delta 
\end{equation*}%
which means that the global error satisfies the tolerance $\delta .$

However, as discussed previously, we do not have $y_{i}$ at each node.
Rather, in the case of local extrapolation, where we have a higher-order
solution available, and we propagate this higher-order solution, we have,
from (\ref{del z z+1}),%
\begin{eqnarray*}
\Delta _{i+1}^{z,z+1} &=&\beta _{i+1}^{z}h^{z+1}+\alpha _{i}^{z,z+1}\Delta
_{i}^{z+1} \\
&=&\varepsilon _{i+1}^{z}+\left( 1+hF_{y}^{z,z+1}\right) \left( \varepsilon
_{i}^{z+1}+\alpha _{i-1}^{z,z+1}\varepsilon _{i-1}^{z+1}+\ldots +\alpha
_{i-1}^{z,z+1}\alpha _{i-2}^{z,z+1}\cdots \alpha _{1}^{z,z+1}\varepsilon
_{1}^{z+1}\right)  \\
&=&\varepsilon _{i+1}^{z}+\sum\limits_{j=1}^{i}\varepsilon
_{j}^{z+1}+O\left( h^{z+3}\right)  \\
&\lesssim &\left\vert \varepsilon _{i+1}^{z}\right\vert
+\sum\limits_{j=1}^{i}\sigma ^{z+2}\left\vert \beta _{j}^{z+1}\right\vert
h^{z+2}=\sigma ^{z+1}\left\vert \beta _{i+1}^{z}\right\vert h^{z+1}+\sigma
^{z+2}i\left\vert \bar{\beta}^{z+1}\right\vert h^{z+2} \\
&=&\sigma ^{z+1}\left\vert \beta _{i+1}^{z}\right\vert h^{z+1}+\sigma
^{z+2}\left\vert \bar{\beta}^{z+1}\left( x_{i}-x_{0}\right) \right\vert
h^{z+1},
\end{eqnarray*}%
where $h$ is assumed to be the stepsize determined from (\ref{h*=...}) and
we have included the safety factor $\sigma $ explicitly. For ease of
analysis, we assume here a uniform stepsize $h$. In the second last line, we
assume that, since $\left\vert \varepsilon _{i+1}^{z}\right\vert \leqslant
\delta $ and $\left\vert \varepsilon _{j}^{z+1}\right\vert \ll \left\vert
\varepsilon _{j}^{z}\right\vert $ (the fundamental assumption in local
extrapolation), we must have $\left\vert \varepsilon _{j}^{z+1}\right\vert
\leqslant \delta .$ We denote the average value of $\beta _{i+1}^{z}$ on $%
\left[ x_{0},x_{i}\right] $ by $\bar{\beta}^{z+1},$ and we use $%
ih=x_{i}-x_{0}$. Now, assuming $\left\vert \bar{\beta}^{z+1}\left(
x_{i}-x_{0}\right) \right\vert h^{z+1}<\delta $ (see the Appendix), we have 
\begin{eqnarray}
\left\vert \Delta _{i+1}^{z,z+1}\right\vert  &\approx &\sigma ^{z+1}\delta
+\sigma ^{z+2}\left\vert \bar{\beta}^{z+1}\left( x_{i}-x_{0}\right)
\right\vert h^{z+1}  \notag \\
&\leqslant &\sigma ^{z+1}\delta +\sigma ^{z+2}\delta   \notag \\
&=&\left( \sigma ^{z+1}+\sigma ^{z+2}\right) \delta .  \label{sigmas}
\end{eqnarray}%
It is easily confirmed that for $\sigma =0.8,$ we have$\left( \sigma
^{z+1}+\sigma ^{z+2}\right) <1$ for $z\geqslant 2.$ This means that%
\begin{equation*}
\left\vert \Delta _{i+1}^{z,z+1}\right\vert <\delta 
\end{equation*}%
for these values of $\sigma $ and $z,$ which suggests that the global error,
like the local error, satisfies the user-defined tolerance. In other words,
propagation of the higher-order solution in local error control via local
extrapolation, has resulted in control of the global error, although the
significance of the safety factor in deriving this result should be clear.
Most importantly, the above result holds only if the assumptions made here
are true; if they are not, then $\left\vert \Delta _{i+1}^{z,z+1}\right\vert 
$ is probably greater than $\delta .$ For this reason, we say that the
global error is \textit{possibly} bounded, but this is \textit{not guaranteed%
}. We should appreciate that such a bounding of the global error is a
beneficial \textit{by-product} of local error control, and is not the
designated objective.

The most important assumption made above is%
\begin{equation*}
\left\vert \bar{\beta}^{z+1}\left( x_{i}-x_{0}\right) \right\vert
h^{z+1}<\delta .
\end{equation*}%
In the Appendix, we show that this is, in fact, a consequence of the
condition%
\begin{equation}
\left\vert \beta _{i+1}^{z}\right\vert h^{z+1}>i\left\vert \bar{\beta}%
^{z+1}\right\vert h^{z+2}.  \label{condition}
\end{equation}%
If this condition is violated, then the assumption does not hold. It is
worth examining this condition in closer detail: the factor $i$ represents
iteration number. It is quite reasonable to assume that there exists a value
of $i$ such that, for the given values of $\bar{\beta}^{z+1}$ and $h,$ we
will have%
\begin{equation*}
\left\vert \beta _{i+1}^{z}\right\vert h^{z+1}<i\left\vert \bar{\beta}%
^{z+1}\right\vert h^{z+2}.
\end{equation*}%
In other words, eventually the RK iterative process will cause (\ref%
{condition}) to be violated. This means that sooner or later the global
error will exceed the imposed tolerance, despite local error control, and
that the bounding of the global error will be `short-lived', so to speak.

We also note that our assumption of a uniform stepsize is reasonable if the
stepsize does not vary considerably; nevertheless, in reality it may do so,
which would also compromise the validity of (\ref{condition}).

\section{Comments}

Some comments are appropriate:

\begin{enumerate}
\item We are not restricted to using a method of order $z+1$ as the
higher-order method in local extrapolation. Any method of order $%
z+r,r\geqslant 1$ can be used. Our analysis and results are essentially
unchanged, save for (\ref{sigmas}), which becomes%
\begin{equation*}
\left\vert \Delta _{i+1}^{z,z+r}\right\vert \leqslant \left( \sigma
^{z+1}+\sigma ^{z+r+1}\right) \delta .
\end{equation*}%
The value of $r$ here will influence the value of the safety factor $\sigma $
for which the coefficient is less than unity.

\item Our analysis clearly shows that global error control cannot be
achieved through local extrapolation. Global error control is usually
achieved through reintegration - estimating the global error after local
extrapolation, and then redoing the computation on the entire interval of
integration with a smaller stepsize. This approach, while effective, can be
inefficient and probably cannot be implemented for real-time problems, where
a globally accurate result is needed before the next iteration. It is not
our intention to report on methods which address this issue, but we would
like to take the opportunity to refer to our own recent work in this regard,
where we have developed an algorithm based on high-order \textit{quenching}
to enable stepwise global error control \cite{jscp2}\cite{jscp3}.
\end{enumerate}

\section{Numerical example}

An instructive example is the initial-value problem%
\begin{eqnarray*}
y^{\prime } &=&\left( \frac{\ln 1000}{100}\right) y \\
y\left( 0\right) &=&1
\end{eqnarray*}%
on $\left[ 0,100\right] .$ The coefficient in the differential equation has
been chosen so that $y$ does not exceed $1000$ on $\left[ 0,100\right] $,
i.e. $y$ does not vary substantially, so that absolute error control is
suitable. The exact solution is $y\left( x\right) =e^{\frac{\ln 1000}{100}%
x}. $

We use RK3 \cite{KC} and RK4 \cite{LV} to implement local error control,
with $\delta =10^{-8}.$ The quantities $\left\vert \varepsilon
^{3}\right\vert $ and $\left\vert \alpha ^{3,4}\Delta ^{4}\right\vert $ are
shown as functions of $x$ in Figure 1 (figure follows appendix). We see that
the local error $\left\vert \varepsilon ^{3}\right\vert $ is bounded by $%
\delta ,$ as expected. However, $\left\vert \alpha ^{3,4}\Delta
^{4}\right\vert $ increases monotonically. The global error $\left\vert
\Delta ^{3,4}\right\vert $ is given by $\left\vert \varepsilon ^{3}+\alpha
^{3,4}\Delta ^{4}\right\vert $ and, although not shown, it is easy to see
that $\left\vert \Delta ^{3,4}\right\vert $ must increase beyond $\delta ,$
and is almost $100$ times greater than $\delta $ at $x=100$. In fact, $%
\left\vert \Delta ^{3,4}\right\vert $ becomes larger than $\delta $ at $%
x=23.6,$ after $118$ iterations. This is the value of $i$ for which (\ref%
{condition}) is violated for this example. For $x<23.6,$ the condition is
not violated and the global error is bounded by $\delta .$ Clearly, though,
this state of affairs does not last, and the propagation of $\Delta ^{4}$
eventually compromises the global accuracy of the solution.

\section{Conclusion}

We have investigated the relationship between local error and global error
in RK methods, under the implementation of local error control via local
extrapolation. We find that, even though local error is successfully
controlled, we cannot expect the same for global error. Our analysis shows
that there is a possibility that the global error will be bounded by a
user-imposed tolerance during the initial stage of the integration, but this
will not last. The propagation of the global error in the higher-order
method will eventually cause the global error in the lower-order method to
exceed the tolerance. A numerical example clearly illustrates these points.

\section{Appendix}

\subsection{The assumption $\left\vert \bar{\protect\beta}^{z+1}\left(
x_{i}-x_{0}\right) \right\vert h^{z+1}<\protect\delta $}

Assuming a constant stepsize $h$, the underlying premise of local
extrapolation is the assumption that 
\begin{equation*}
\left\vert \beta _{i+1}^{z}\right\vert h^{z+1}\gg \left\vert \beta
_{i+1}^{z+1}\right\vert h^{z+2},
\end{equation*}%
which implies%
\begin{equation*}
\left\vert \beta _{i+1}^{z}\right\vert h^{z+1}=M\left\vert \beta
_{i+1}^{z+1}\right\vert h^{z+2},
\end{equation*}%
where $M\gg 1$ is a large number. Assuming that $\beta _{i+1}^{z+1}$ is a
slowly varying function of $x,$ we can replace $\beta _{i+1}^{z+1}$ with its
average value $\bar{\beta}^{z+1},$ and so we may write%
\begin{equation*}
\left\vert \beta _{i+1}^{z}\right\vert h^{z+1}>i\left\vert \bar{\beta}%
^{z+1}\right\vert h^{z+2}=\left\vert \bar{\beta}^{z+1}\left(
x_{i}-x_{0}\right) \right\vert h^{z+1},
\end{equation*}%
subject to the condition that $i<M$. Hence,%
\begin{equation*}
\left\vert \beta _{i+1}^{z}\right\vert h^{z+1}<\delta \Rightarrow \left\vert 
\bar{\beta}^{z+1}\left( x_{i}-x_{0}\right) \right\vert h^{z+1}<\delta .
\end{equation*}

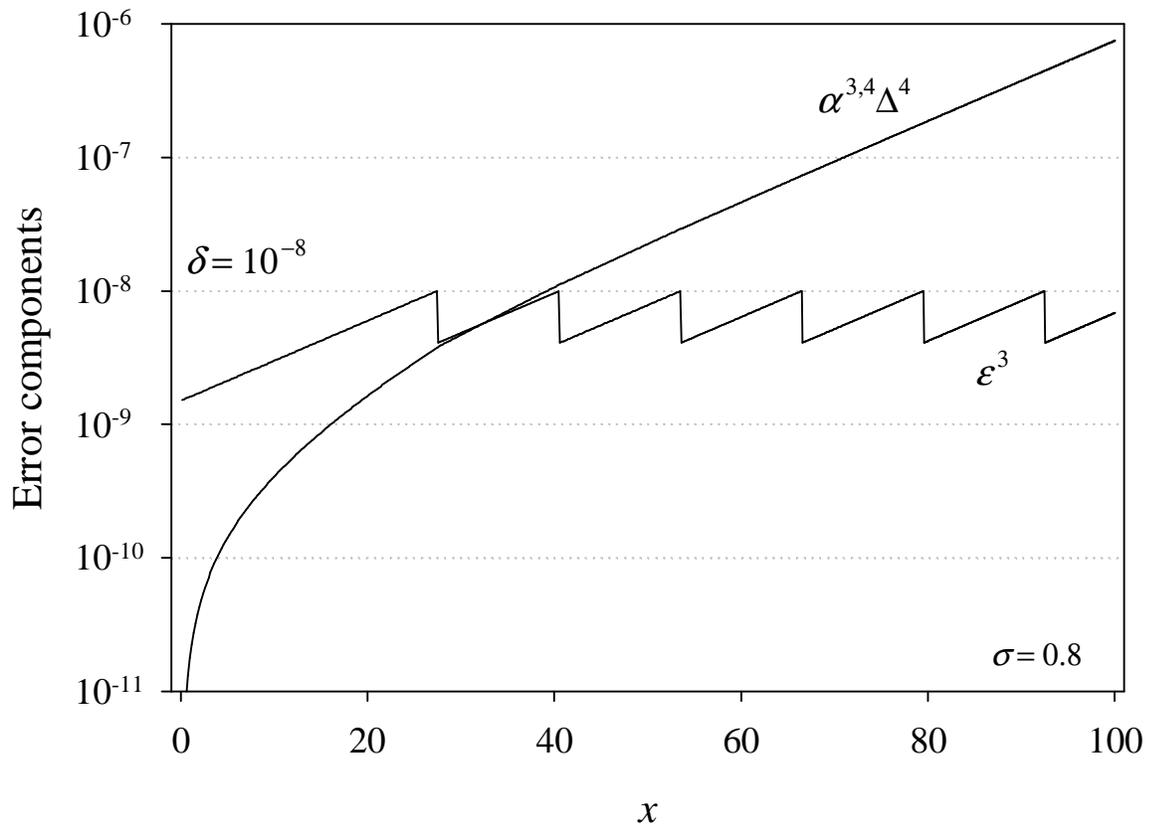
\end{document}